\documentclass{amsart}
\usepackage{amssymb}
\newtheorem{theorem}{Theorem}
\title{Beyond Uncountable}
\author{Paola Cattabriga}
\address{University of Bologna, Italy}
\email{co14099@iperbole.bologna.it}
\begin{document}
\maketitle
\thispagestyle{empty}
\bigskip \bigskip
\hspace{1.84in} \parbox{2.90in}{\scriptsize 
$\dots$
The fact is that such a procedure is not   applicable. Why? Because their definitions
are not predicative and  contain within such a vicious circle  I already mentioned  
above; not predicative definitions can not be substituted to 
defined terms. In this condition,
\emph{logistics is no longer sterile: it generates contradictions}.
(Jules-Henri Poincar\'e 1902, \cite{poincare} 211, our translation.)}
\bigskip \bigskip
\section*{Introduction}
By common consent  Russell's antinomy is the reason 
why in  Zermelo--Fraenkel set theory,  
there is no set which comprehends \emph{all} sets. Furthermore, given
any set $A$, there is no set which contains \emph{all} sets  which are not 
members of $A$ (in particular, there is no set which is the complement of $A$) 
(\cite{foundations} 40-41). In other words, given any set $A$, the \emph{absolute}
complement of $A$, i.e. $\{x\mid x\notin A\}$, cannot be defined and the complement
of $A$, can only be defined  as  \emph{relative} to another given set.
For instance, if $A$ is a subset of $B$, then the relative complement of $A$ in $B$ is 
defined by $$B - A =\{x \in B\mid x\notin A\}.$$ The existence of the relative 
complement is ensured by the axiom schema of the Subsets
\begin{equation}\label{sub}
\forall z_{1}\dots z_{n}\forall s \exists  y \forall x
(x\in y \iff x \in s \wedge\varphi(x)),
\end{equation}
where $\varphi(x)$ is a first order well formed formula, $ z_{1},\dotsc,z_{n},x$ are
the free variables of $\varphi(x)$, and $y$ is not free in $\varphi(x)$, which admits
general comprehension only for members $x$ of a given set $s$. Indeed we
are always allowed to assert
\begin{equation}\label{csub}
\forall z\forall s \exists  y \forall x
(x\in y \iff x \in s \wedge x \notin z),
\end{equation}
as an instance of (\ref{sub}). This set $y$ is the relative complement of $z$ in $s$
(\cite{fraenkel} 23). This premise and the following subsection are 
introductory to the results of Section \ref{s1}. 
In this abstract Zermelo--Fraenkel set theory stands for general 
 first order set theory. 

\subsubsection*{Basic setup} We refer to  Zermelo--Fraenkel set 
theory   with  $\mathfrak{ZF}$.
Let us  recall the axiom of Extensionality  
\begin{equation}\label{ext}
\forall x\forall y\:[
\forall z(z\in x \iff z \in y) \Longrightarrow x = y)],
\end{equation}
and the other main concepts we  shall be concerned with. For  details see
\cite{fraenkel,foundations}.
$$x \subseteq z\: = \:\forall w\:(w\in x \Longrightarrow w\in z),$$
$${\mathit P}(z) \:= \:\{x \mid x\subseteq z\},$$
$x \sim y$ denotes that the sets $x$ and $y$ are equinumerous or equal in cardinality, 
namely there exists a one to one correspondence between their elements.

\medskip
\noindent $x \not\sim y$ denotes that the sets $x$ and $y$ are not equal in cardinality, 
namely there exists no one to one correspondence between their elements.

\medskip
\noindent $x <_{c} y$ and $x \leq _{c} y$ denote respectively that the set $x$ has
cardinality properly less than the cardinality of $y$, and that the set
$x$ has cardinality less than $y$. 

\medskip
Let us also recall the argumentation of the so-called  Cantor's  theorem. 
We shall present the version in (\cite{moschovakis} 15), 
for a more detailed exposition the reader is referred to 
(\cite{cantor1, cantor2},\cite{fraenkel}). 
\subsubsection*{(Cantor's proposition)} 
{\it For every set A in $\mathfrak{ZF}$, $$ A <_{c}{\mathit P}(A)$$ i.e. $A \leq _{c} {\mathit P}(A)$ but 
$A \not\sim  {\mathit P}(A)$.} 

\medskip

That $A \leq _{c} {\mathit P}(A)$ follows from the fact that the function
$$x \mapsto \{x\}$$ which associates with each member $x$ of $A$ its singleton 
$\{x\}$ is an injection. To complete the proof, we assume, toward
a contradiction, that there exists a one to one correspondence
$$\mathtt{g}:A \mapsto {\mathit P}(A)$$
which establishes that $A \sim {\mathit P}(A)$ and we define the set
\begin{equation}\label{part1}
B \; = \; \bigl\{ x \in A \mid x \notin \mathtt{g}(x)\bigr\}.
\end{equation}
Now $B$\footnote{Notice that in $\mathfrak{ZF}$ the definition of 
$B$ is an example of (\ref{sub}), as 
one can easily see, $B$ is defined within $A$.} is a subset of $A$ and $\mathtt{g}$ 
is a surjection, so \emph{there
must exist some $b\in A$ such that $B = \mathtt{g}(b)$} ({\em diagonalization}), 
and (as for each
$x\in A$) either  $b\in B$ or $b\notin B$.
\begin{itemize}
\item[(*)] If $b\in B$ then $b\in \mathtt{g}(b)$ since $B = \mathtt{g}(b)$,
so that $b$ does not satisfy the condition which defines $B$, and hence
$b\notin B$, contrary to hypothesis.
\item[(**)] If $b\notin B$, then $b\notin \mathtt{g}(b)$, so that $b$ now satisfies
the defining condition for $B$ and hence $b\in B$, which again contradicts the hypothesis.
\end{itemize} 
Thus we reach a contradiction from the assumption that the bijection $\mathtt{g}$
exists and the proof is complete.  

\section{The relative complement}\label{s1}
Let us read the above so-called Cantor's theorem
and connect again to (\ref{part1}), i.e. the step of the definition of $B$ within Cantor's  argumentation. 
As previously observed,  the relative complement can always be defined,
thanks to (\ref{csub}). Accordingly let us define $ \overline{B}\; =\; A - B$
as  the relative complement of $B$ in $A$, i. e.
\begin{equation}\label{part1-2}
\overline{B} \; = \; \{ x \in A \; | \; x \not\in B\}. 
\end{equation}
One can easily see that 
\begin{equation}\label{part2}
 \overline{B} \; = \; \{ x \in A \; | \; x \in \mathtt{g}(x)\}. 
\end{equation}
In other words, being (\ref{part1-2}) legitimated by (\ref{csub}), 
whenever (\ref{part1}) is defined immediately (\ref{part2}) is  
defined too.

Consequently we have in $\mathfrak{ZF}$ the following situation
\begin{xxalignat}{3}
\qquad\quad&\notag\mathtt{g}:A \mapsto {\mathit P}(A) & \qquad\qquad\qquad\qquad\text{by assumption,}\\
\qquad\quad&\notag B \; = \; \bigl\{ x \in A \mid x \notin \mathtt{g}(x)\bigr\} &
 \qquad\qquad\qquad\qquad\text{by Subset axiom,}\\
\qquad\quad&\notag \overline{B} \; = \; \bigl\{ x \in A \mid x \in \mathtt{g}(x)\bigr\} &
\qquad\qquad\qquad\qquad\text{by Subset axiom.}
\end{xxalignat}
We can then state that $B$ and $\overline{B}$ are subsets of $A$. 
By its 
definition $\mathtt{g}$ is a surjection and for each $x\in A$ we have either $x\in B$ or
$x\notin B$, i.e. by (\ref{part1-2}) either $x\in B$ or $x\in \overline{B}$.
Let us reconsider the statement
 \emph{there must exist some $b\in A$ such that $B = \mathtt{g}(b)$}
 ({\em diagonalization}), within  Cantor's  argumentation.
If such $b$ exists, from $B\neq \overline{B}$, we obtain $B = \mathtt{g}(b) \leftrightarrow 
\overline{B} \neq \mathtt{g}(b)$, i.e. $B = \mathtt{g}(b) $ or 
$\overline{B} = \mathtt{g}(b)$ but not both.
 We  have then the main consequence of taking into consideration the 
definition of the relative 
complement with respect to Cantor's argumentation in $\mathfrak{ZF}$. 
Applying (\ref{ext}) we obtain
\begin{equation}
(b\in \overline{B} \iff b \in \mathtt{g}(b)) \Longrightarrow 
\overline{B} = \mathtt{g}(b),
\end{equation}
 hence  by (\ref{part2})
\medskip
\begin{equation}
\overline{B} = \mathtt{g}(b).
\end{equation}
\medskip
Moreover since $b\in B$ or $b\in \overline{B}$ but not both, and $B=\mathtt{g}(b)$ or 
$\overline{B}=\mathtt{g}(b)$ but not both we have
\medskip
\begin{equation}
B \neq \mathtt{g}(b).
\end{equation}
\medskip
Accordingly the assertion \emph{there must exist some $b\in A$ such that 
$B =\mathtt{g}(b)$} is false. By the axiom schema of Subsets and the axiom of Extensionality, 
diagonalization can not be stated as true in $\mathfrak{ZF}$. 
Consequently (*) and (**) cannot be accomplished  and Cantor's theorem does not hold 
in $\mathfrak{ZF}$. 
In fact we have only two cases

\bigskip

\begin{tabular}{lp{3in}}
1. $b\in B$ and $\overline{B} = \mathtt{g}(b)$, &then $b\notin \mathtt{g}(b)$ so that
$b$ satisfies the condition in (\ref{part1}) which defines $B$, and hence $b\in B$,
accordingly to the hypothesis;\\
2. $b\in \overline{B}$ and $\overline{B} = \mathtt{g}(b)$, &then $b \in \mathtt{g}(b)$
so that $b$ satisfies condition in (\ref{part2}), and hence $b \in \overline{B}$,
accordingly to the hypothesis.
\end{tabular}

\medskip

\noindent We have thus established  the following theorem.

 \begin{theorem}
By the definability of the relative complement, 
 Cantor's proposition does not hold as a theorem in $\mathfrak{ZF}$.
\end{theorem}

\section{The restriction on uniqueness}\label{s2}

Let us leave aside now Cantor's argumentation. 
We assume to have a set $A$ already defined in $\mathfrak{ZF}$. By 
the axiom schema of the Subsets we have
\begin{equation}\tag{I}
\vdash_\mathfrak{ZF} (b\in B \Longleftrightarrow b\notin 
\mathtt{g}(b) \wedge b\in A),
\end{equation}
which defines $B$ as a subset of $A$.
Since $b \in A$ is true we obtain
\begin{equation}\tag{II}\vdash_\mathfrak{ZF} (b\in B \Longleftrightarrow b\notin \mathtt{g}(b)).
    \end{equation}
Furthermore by the axiom of Extensionality and the underlying laws for 
identity  ($\forall z (x\in z \Longleftrightarrow y\in z) \Longleftrightarrow 
x = y$,   \cite{foundations} 25, 28)

\begin{equation}\tag{III}
\vdash_\mathfrak{ZF} (b\in B \iff b \in \mathtt{g}(b)) \iff B = \mathtt{g}(b),
\end{equation}
and therefore 
\begin{equation}\tag{IV}
\vdash_\mathfrak{ZF} (b\in B \iff b \notin \mathtt{g}(b)) \Longrightarrow 
B \neq \mathtt{g}(b),
\end{equation}
so that by (II) and (IV) 
\medskip
\begin{equation}\tag{V}
\vdash_\mathfrak{ZF}  B \neq \mathtt{g}(b).
\end{equation}

Let us   state in $\mathfrak{ZF}$
\begin{equation}\tag{VI} 
    B = \mathtt{g}(b),
\end{equation}
then we  attain
\begin{equation}\tag{VII}  
    \vdash_\mathfrak{ZF} B = \mathtt{g}(b)\:  \wedge  \:B \neq \mathtt{g}(b),
    \end{equation}
 accordingly $\mathfrak{ZF}$ turns out to be inconsistent. In simple 
 terms, to
state diagonalization, $B = \mathtt{g}(b)$, as true makes $\mathfrak{ZF}$ inconsistent. 
There is no need to yield diagonalization within the contest of a reasoning 
or argumentation.
A definition like (I) leads  to contradiction in any case.  The 
explanation can be provided by the theory of definition which states 
the conditions and restrictions for defining proper equivalence in mathematics 
(see for example \cite{suppes} 151-173).
Definition (I)  neglects a restriction embodied in the rules for proper 
definitions, established on the basis of
the criterions of eliminability and non-creativity. Exactly as 
Russell's antinomy, definition of $B$ drops the restriction on 
uniqueness,  
which is given when defining a new operation symbol (or a new individual 
constant, i.e. an operation symbol of rank zero) \cite{pc2}. An 
equivalence like $$O(x_{1}\ldots x_{n}) = y \iff \Phi$$ introducing a new 
operation symbol $O$, is a proper definition only if the formula 
$$\exists ! y \Phi$$ is derivable from the axioms and preceding 
definitions of the theory (\cite{suppes} 158-159). In 
$\mathfrak{ZF}$ the uniqueness is ensured by the axiom of 
Extensionality (\ref{ext}), which implies that there exists \emph{at 
most}, one set $y$, which contains exactly those elements $x$ which 
fulfill the condition $\varphi (x)$ in (\ref{sub}) (\cite{foundations} 31).
If $B$ and $\mathtt{g}(b)$ are two sets each of which contains 
exactly those elements $b$ which fulfil the condition $b\in 
\mathtt{g}(b)$, then $B$ and $\mathtt{g}(b)$ are equal, see (III). 
Accordingly, there exists \emph{at most}, one $B$, such that $b\in 
\mathtt{g}(b)$. Definition (I) implying $B \neq \mathtt{g}(b)$, (IV) 
and (V), neglects the restriction on uniqueness established by 
Extensionality and therefore the relative consistency
   embodied in the criterion of non-creativity  (\cite{suppes} 155; \cite{pc2}). 
This explains why Extensionality blocks the derivation of the 
existence of some $b\in A$ such that $B = \mathtt{g}(b)$ in Cantor's 
argumentation. Moreover, this fulfils the criterion established by an 
editor, according to which, \emph{to attack an argument, you must find 
something wrong in it}. We showed indeed that the definition of $B$, 
neglecting the restriction on uniqueness, is 
always wrong in $\mathfrak{ZF}$ and therefore a wrong \emph{object sentence}  
in Cantor's argumentation \cite{hodges}. 

When
 this results  are regarded  together with those presented in 
 \cite{pc1, pc3} it arises clearly a  similitude.  If  a set, 
 or a predicate, is  object of diagonalization then the definition
of its complement  leads  to the  invalidity of
the diagonalization itself. In Section \ref{s3} we shall apply this code of behavior to 
Cantor's diagonal argument.

\section{Cantor's diagonal argument}\label{s3}
In 1891 Cantor presented a striking argument which has come to know 
as Cantor's diagonal argument  \cite{cantor2}. It  
runs as follows.

Consider the two elements $m$ and $v$.
Let $M$ be the set whose elements $E$ are sequences 
$< x_{1},x_{2},\ldots , x_{v} , \ldots  >$
where each of $ x_{1},x_{2},\ldots , x_{v} ,$ $ \ldots $ is either $m$ 
or $w$.

{\bf Cantor's proposition } 
{\it If $E_{1},$ $E_{2},$ $\ldots , E_{v} , \ldots $ is any simply infinite 
sequence of elements of the set $M$, then there is always an element 
$E_{0}$ of $M$ which corresponds to no $E_{v}$.}

To prove this proposition, Cantor arranged a denumerable list of 
elements in an array.

\bigskip

\begin{align}
 \notag   E_{1} &  =    < a_{1,1}, a_{1,2} \ldots , a_{1,v} , \ldots >  \\
 \notag	E_{2} &  =   < a_{2,1}, a_{2,2} \ldots , a_{2,v} , \ldots >\\
 \notag	  \quad &. \quad . \quad . \quad . \quad . \quad . \quad . \quad . \quad .  \\
 \notag	E &  =   < a_{\mu,1}, a_{\mu,2} \ldots , a_{\mu,v} , \ldots >\\
 \notag	  \quad &. \quad . \quad . \quad . \quad . \quad . \quad . \quad . \quad .  \\ 
\notag	\end{align}
Each $a_{\mu,v}$ is either $m$ or $w$. Cantor defined a sequence 
$ b_{1},b_{2},b_{3},\ldots , $ where each element is $m$ or $w$, and,  if 
$a_{v,v} = m$ then $b_{v} = w$, and if $a_{v,v} = w$ then $b_{v} = m$.
Let $E_{0}   =    < b_{1}, b_{2} , b_{3} , \ldots >$. Then no $E_{v}$ 
corresponds to $E_{0}$, by reason that $b_{v}\neq a_{v,v}$.  

\bigskip 

$E_{1},$ $E_{2},$ $\ldots , E_{v} , \ldots $ is {\it any simply infinite sequence of 
elements of the set $M$}, so that we can think to  a definite infinite sequence 
of elements of $M$ as follows.

\bigskip

\begin{align}
 \notag   E_{1}^{*} &  =    < a_{1,1}^{*}, a_{1,2}^{*} \ldots , a_{1,v}^{*} , \ldots >  \\
 \notag	E_{2}^{*} &  =   < a_{2,1}^{*}, a_{2,2}^{*} \ldots , a_{2,v}^{*} , \ldots >\\
 \notag	  \quad &. \quad . \quad . \quad . \quad . \quad . \quad . \quad . \quad .  \\
 \notag	E^{*} &  =   < a_{\mu,1}^{*}, a_{\mu,2}^{*} \ldots , a_{\mu,v}^{*} , \ldots >\\
 \notag	  \quad &. \quad . \quad . \quad . \quad . \quad . \quad . \quad . \quad .  \\ 
\notag	\end{align}

\bigskip

Each $a_{\mu,v}^{*}$ is either $m$ or $w$ and if $a_{\mu,v}= 
m$ then $a_{\mu,v}^{*} =w$, if $a_{\mu,v} = w$ then $a_{\mu,v}^{*}=m$.
Then $b_{v} = a_{v,v}^{*}$ and
$E_{0}$ is never different on the $v$-th coordinate, so that it could 
even be for some $E_{v}^{*}$ that
$$E_{v}^{*} = E_{0}.$$
Since $E_{1}^{*}, E_{2}^{*}, \ldots E^{*} \ldots$ is a simply 
infinite sequence of elements of $M$, previous Cantor's proposition is false.

\bigskip 

In both the cases of Cantor's power set theorem and Cantor's diagonal 
argument,  the definition  of the complement of the 
object of diagonalization 
leads to the rejection of the diagonalization itself.

Working on logical complementation, Section \ref{ac}, 
gives  proof of the Axiom of Choice in $\mathfrak{ZF}$, and 
its
refutation in a framework which is  no longer $\mathfrak{ZF}$,  on the basis of the 
universal validity of the first order logical 
truths.

\section{The Axiom of Choice}\label{ac}

Let us consider the following first order logic formula
\begin{equation}\label{newp}
\forall z\:\forall y\: \forall x\:\bigl[ ( x \in y \iff
x \subseteq z \bigr)  \: \iff\:(x \notin y \iff
\neg(x \subseteq z))\bigr].
\end{equation}
One can easily see it is a logically valid formula.
We can then assume (\ref{newp})  holds in $\mathfrak{ZF}$ (Zermelo Fraenkel 
set theory, \cite{foundations}).
From a comparison with the classical  Axiom of Power Set 
\begin{equation}\label{power}
\forall z\:\exists y\: \forall x\:\bigl( x \in y \iff
x \subseteq z \bigr)
\end{equation}
where $y$ is ${\mathit P}(z)$, the power-set of $z$,
it follows immediately  that  (\ref{newp}) holds for
the power-set because of its logical validity. We can then assume that (\ref{newp})  
 establishes the definition of ${\mathit P}(z)$ as inseparable from the 
definition of its complement. We  could think of $\overline{{\mathit P}}(z)$    
as defined by 
\begin{equation}\label{nopower}
\forall z\:\exists y\: \forall x\:\bigl( x \in y \iff
\neg (x \subseteq z) \bigr).
\end{equation}

Let us recall the Axiom of Choice as defined in (\cite{foundations} 39-40, 53-55),
to which the reader is referred for  details. 
\medskip
\begin{itemize}
\item[(AC)] \emph{If $t$ is a disjointed set which does not contain the null--set, 
its outer product $\mathcal{OP}t$ is different from the null--set}.
\end{itemize}
\medskip
In other words, among the subsets of $\bigcup t$ there is at least one whose
intersection with each member of $t$ is a singleton.

$\mathcal{OP}t$ exists only if $\mathit{P}(\bigcup t)$, the set of the subsets of 
the union of $t$, exists. Immediately, by (\ref{newp}), AC is true, since, if 
$t$ does not contain the null--set, $\mathit{P}(\bigcup t)$ is never disjointed, and
therefore there are selection sets of $t$ and $\mathcal{OP}t$ is different from 
the null--set.

The reasons lie in the logical structure of (\ref{newp}), which states 
$\mathit{P}(z)$, namely $y$,  to be a set excluding those parts $x$ for which  
$\neg\:(x\subseteq z)$ holds (see the component
$(x \notin y \iff \neg\:(x\subseteq z))$ in (\ref{newp})). 

Let us consider (\ref{newp}), with $\emptyset\notin z$. $\neg\:(x\subseteq z)$ is 
true if $(z \subset x) \:\vee\:(x\cap z=\emptyset)$, hence $\mathit{P}(z)$ is constructed
in such a way that $x$ and $z$ are never disjointed (since
$(\neg (z \subset x) \:\wedge\:\neg(x\cap z=\emptyset))$ holds for  $\mathit{P}(z)$ 
by (\ref{newp})).
 
Now, consider $t$ in AC. $\emptyset\notin t$, hence all the members of  
$\mathit{P}(\bigcup t)$ are never disjointed, also if $t$ is disjointed. 
Consequently $\mathcal{OP}t$ is different from the null--set.
 
The  opposite holds for $\overline{\mathit P}$, because the members of 
$\overline{\mathit{P}}(\bigcup t)$ are always disjointed, and even if $t$ does not 
contain the null--set, its outer product is equal to the empty set. To visualize how 
$\overline{\mathit{P}}$  gives rise to this situation, we can consider the following 
logically valid formula
\begin{equation}\label{negnewp}
\forall z\:\forall y\: \forall x\:\bigl[ ( x \in y \iff
\neg(x \subseteq z) )  \: \iff\:(x \notin y \iff
\neg\neg(x \subseteq z))\bigr].
\end{equation}
Since $\emptyset\notin z$, we can think of $y$ as excluding those
parts $x$ such that $x\subseteq z$ 
(and $\neg ((z \subset x) \:\vee\:(x\cap z=\emptyset))$). Thus for the parts
$x$ of $y$ it holds always $(x\cap z=\emptyset)$.

To summarize we have  proved that
\medskip
\begin{equation}\tag*{}
 [\text{(1) and }\mathit{P}(z)] \,\Longrightarrow\, \text{AC} \qquad
 \qquad
[ \text{(4) and }\overline{\mathit P}(z) ]\,\Longrightarrow \,\text{not AC.}
\end{equation}
\medskip
We can then assert  the following theorems.
\begin{theorem}
The Axiom of Choice holds by (\ref{newp}). 
\end{theorem}
\begin{theorem}
The negation of the Axiom of Choice holds   by (\ref{negnewp}) 
(if  $\overline{\mathit{P}}$  is defined).
\end{theorem}

{\Tiny \it Copyright \copyright 1998-2006 Paola Cattabriga}
\end{document}